\documentclass[12pt]{article}
\usepackage{amssymb}
\usepackage{amsfonts}
\usepackage{times}
\usepackage{mathptmx}
\usepackage{amsmath}
\usepackage[usenames]{color}
\usepackage{mathrsfs}
\usepackage{amsfonts}
\usepackage{amssymb,amsmath}
\usepackage{CJK}
\usepackage{cite}
\usepackage{cases}
\usepackage{amsthm}

\def\dis{\displaystyle}\def\ssc{\scriptscriptstyle}\def\OP#1{\raisebox{-7pt}{$\stackrel{\dis\oplus}{\ssc #1}$}}
\def\deg{{\rm deg}}

\def\cl{\centerline}

\def\a{\alpha}

\def\b{\beta}
\def\vs{\vspace*}

\def\Z{\mathbb{Z}}

\def\C{\mathbb{C}}

\def\QED{\hfill$\Box$}
\def\ni{\noindent}

\numberwithin{equation}{section}
\newtheorem{theo}{Theorem}[section]
\newtheorem{defi}[theo]{Definition}

\newtheorem{lemm}[theo]{Lemma}
\newtheorem{prop}[theo]{Proposition}

\begin{document}
\begin{center}
{\bf\large Loop Heisenberg-Virasoro Lie Conformal algebra}
\footnote {Supported by NSF grant no.~11431010, 11371278, 11001200 and 11101269 of China, the Fundamental Research Funds for the Central Universities, the grant no.~12XD1405000 of Shanghai Municipal Science and Technology Commission.

Corresponding author:  Y.~Su (ycsu@tongji.edu.cn)
}
\end{center}

\cl{Guangzhe Fan$^{\,*}$, Yucai Su$^{\,*}$,
Henan Wu$^{\,\dag}$}

\cl{\small $^{\,*}$Department of Mathematics, Tongji University, Shanghai 200092, China}
\cl{\small $^{\,\dag}$ School of Mathematical Sciences, Shanxi University, Taiyuan 030006, China}

\vs{8pt}

{\small\footnotesize
\parskip .005 truein
\baselineskip 3pt \lineskip 3pt
\noindent{{\bf Abstract:} Let $HV$ be the loop Heisenberg-Virasoro Lie algebra over $\C$ with basis $\{L_{\a,i},H_{\b,j}\,|\,\a,\,\b,i,j\in\Z\}$ and brackets $[L_{\a,i},L_{\b,j}]=(\a-\b)L_{\a+\b,i+j}, [L_{\a,i},H_{\b,j}]=-\b H_{\a+\b,i+j},[H_{\a,i},H_{\b,j}]=0$. In this paper, a formal distribution Lie algebra of $HV$ is constructed. Then the associated  conformal algebra $CHV$ is studied, where $CHV$ has a $\C[\partial]$-basis $\{L_i,H_i\,|\,i\in\Z\}$ with $\lambda$-brackets $[L_i\, {}_\lambda \, L_j]=(\partial+2\lambda) L_{i+j}, [L_i\, {}_\lambda \, H_j]=(\partial+\lambda) H_{i+j}, [H_i\, {}_\lambda \, L_j]=\lambda L_{i+j}$ and $[H_i\, {}_\lambda \, H_j]=0$. In particular, the conformal derivations of $CHV$ are determined. Finally, rank one conformal modules and $\Z$-graded free intermediate series modules over $CHV$ are classified. \vs{5pt}

\parskip .001 truein\baselineskip 6pt \lineskip 6pt

}}
\parskip .001 truein\baselineskip 6pt \lineskip 6pt

\section{Introduction}
Lie conformal algebra encodes an axiomatic description of the operator product expansion of chiral fields in conformal field theory. Kac introduced the notion of the conformal algebra in \cite{K3}. Conformal module is a basic tool for the construction of free field realization of infinite dimensional Lie (super)algebras in conformal field theory. In recent years, the structure theory, representation theory and cohomology theory of Lie conformal algebras have been  extensively studied by many scholars. For example, a finite simple Lie conformal algebra was proved to be isomorphic to either the Virasoro conformal algebra or  the current conformal algebra associated with a finite dimensional simple Lie algebra in \cite{DK}. Finite irreducible conformal modules over the Virasoro conformal algebra were determined in \cite{CK}. The cohomology theory of conformal algebras was developed in  \cite{BKV}. The low dimensional cohomologies of the infinite rank general Lie conformal algebras $gc_N$ with trivial coefficients were computed in \cite{S}. Two new nonsimple conformal algebras associated with the Schr$\ddot{\rm o}$dinger-Virasoro Lie algebra and the extended Schr$\ddot{\rm o}$dinger-Virasoro Lie algebra were constructed in \cite{SY1}. The Lie conformal algebra of a Block type was introduced and free intermediate series modules were classified in \cite{GXY}. The loop Virasoro Lie conformal algebra was studied in \cite{WCY}.

The Heisenberg-Virasoro algebra contains the classical Heisenberg algebra and the Virasoro algebra as subalgebras. As the universal central extension of the Lie algebra of differential operators on a circle of order at most one, the Heisenberg-Virasoro algebra has been widely studied in the mathematical and physical literature. For example, the twisted Heisenberg-Virasoro algebra has been first studied by E. Arbarello et al. in \cite{A}. Various generalizations of the Heisenberg-Virasoro algebra have been extensively studied by several authors (e.g., \cite{LZ,B}). However, it seems to us that little has been known on the loop Heisenberg-Virasoro algebra.

Infinite rank Lie conformal algebras are important ingredients of Lie conformal algebras. In this paper, we 
study an infinite rank Lie conformal algebra, namely, the loop Heisenberg-Virasoro Lie conformal algebra $CHV$ (cf.~\eqref{1.3}). The loop Heisenberg-Virasoro algebra $HV$ is defined to be a Lie algebra with basis $\{L_{\a,i},H_{\b,j}\,|\,\a,\,\b,$ $i,j\in\Z\}$ and Lie brackets given by
\begin{eqnarray*}
\aligned
&[L_{\a,i},L_{\b,j}]=(\a-\b)L_{\a+\b,i+j}, \ \ 
&[L_{\a,i},H_{\b,j}]=-\b H_{\a+\b,i+j}, \ \ \ \ 
&[H_{\a,i},H_{\b,j}]=0.
\endaligned
\end{eqnarray*}
Let $W$ be the subalgebra spanned by
$\{L_{\a,i}\,|\,\a,i\in\Z\}$ and
$H$ the ablian ideal spanned by
$\{H_{\b,j}\,|\,\b,j\in\Z\}$.
Then $W$ is actually isomorphic to the centerless loop Virasoro algebra, and $HV$ also contains the Heisenberg-Virasoro
algebra
\begin{equation*}
\mathscr{HV}={\rm span}\{L_{\a,0},\,H_{\b,0}|\,\a,\,\b\in\Z\}.
\end{equation*}
The Lie conformal algebra of $HV$, denoted by $CHV$, is constructed in Section $2$.
As one can see, it is a Lie conformal algebra with $\C[\partial]$-basis $\{L_i,H_j\,|\,i,\,j\in\Z\}$ and $\lambda$-brackets
\begin{eqnarray}\label{1.3}
\aligned
&[L_i\, {}_\lambda \, L_j]=(\partial+2\lambda) L_{i+j}, \ \
[L_i\, {}_\lambda \, H_j]=(\partial+\lambda) H_{i+j},\ \  [H_i\, {}_\lambda \, L_j]=\lambda H_{i+j}, \ \
[H_i\, {}_\lambda \, H_j]=0.
\endaligned
\end{eqnarray}
We 
remark that the conformal subalgebra $CVir=\C[\partial]L_{0}$ is isomorphic to the well-known Virasoro conformal algebra and
the conformal subalgebra\begin{equation}\label{CWAAAA}CW=\OP{i\in\Z}\C[\partial]L_i,\end{equation} is isomorphic to the loop Virasoro Lie conformal algebra studied in \cite{WCY}.
So some results about this Lie conformal algebra can be applied in this case.

This paper is organized as follows. In Section $2$,  some basic definitions of Lie conformal algebras are  recalled. In Section $3$, we 
start from $HV$ to construct its Lie algebra of formal distributions. Then we 
construct the related Lie conformal algebra $CHV$. In Section $4$, conformal derivations of $CHV$ are determined. Finally, 
rank one conformal modules and $\Z$-graded free intermediate series modules over $CHV$ are classified in Section $5$ and Section $6$.

Throughout the paper, we denote by $\C,\,\C^*,\, \Z$ the sets of complex numbers, nonzero complex numbers, integers respectively.
\section{Preliminaries}
In this section, we recall some definitions related to Lie conformal algebras in \cite{DK,K1,K3}.
A formal distribution (usually called a field by physicists) with coefficients in a complex
vector space $U$ is a series of the following form:
\begin{equation*}
a(z)=\mbox{$\sum\limits_{i\in\Z}$}a_{(i)} z^{-i-1},
\end{equation*}
where $z$ is an indeterminate and $a_{(i)}\in U$.
Denote by $U[[z,z^{-1}]]$ the space of formal distribution with coefficients in $U$.
The space $U[[z,z^{-1},w,w^{-1}]]$ is defined in a similar way. A formal distribution $a(z,w)\in U[[z,z^{-1},w,w^{-1}]]$ is called {\it local} if $(z-w)^N a(z,w)=0$ for some $N\in\Z^+$.
Let $g$ be a Lie algebra.
Two formal distributions $a(z), b(z)\in g[[z,z^{-1}]]$ are
called {\it pairwise local} if $[a(z),b(w)]$ is local in $g[[z,z^{-1},w,w^{-1}]]$.
\begin{defi}\rm
A family $F$ of pairwise local formal distributions, whose coefficients span $g$, is called
a formal distribution Lie algebra of $g$. In such a case, we say that the family $F$ spans $g$. We will write $(g, F)$ to
emphasize the dependence on $F$.
\end{defi}

Define the {\it formal delta distribution} to be
\begin{equation*}
\delta(z,w)=\mbox{$\sum\limits_{i\in\Z}$}z^i w^{-i-1}.
\end{equation*}
The following proposition describes an equivalent condition for a formal
distribution to be local.
\begin{prop}\label{local}
A formal distribution $a(z,w)\in U[[z,z^{-1},w,w^{-1}]]$ is local if and only if $a(z,w)$ can be written as
\begin{equation*}
a(z,w)=\mbox{$\sum\limits_{j\in\Z^+}$}c^j(w)\frac{\partial^j_w\delta(z,w)}{j!}
\mbox{ \ (finite sum) \ for some $c^j(w)\in U[[w,w^{-1}]]$.}\end{equation*}
\end{prop}

In this paper, we adopt the following definition of Lie conformal algebras using $\lambda$-brackets as in\cite{K3}.
\begin{defi}\label{D1}\rm
A Lie conformal algebra is a $\C[\partial]$-module $A$ endowed with a $\lambda$-bracket $[a{}\, _\lambda \, b]$ which defines a linear map $A\otimes A\rightarrow A[\lambda]$, where $\lambda$ is an indeterminate and $A[\lambda]=\C[\lambda]\otimes A$, subject to the following axioms:
\begin{equation}\label{conformal}
\aligned
&[\partial a\,{}_\lambda \,b]=-\lambda[a\,{}_\lambda\, b],\ \ \ \
[a\,{}_\lambda \,\partial b]=(\partial+\lambda)[a\,{}_\lambda\, b];\\
&[a\, {}_\lambda\, b]=-[b\,{}_{-\lambda-\partial}\,a];\\
&[a\,{}_\lambda\,[b\,{}_\mu\, c]]=[[a\,{}_\lambda\, b]\,{}_{\lambda+\mu}\, c]+[b\,{}_\mu\,[a\,{}_\lambda \,c]].
\endaligned
\end{equation}
\end{defi}

For any local formal distribution $a(z,w)$, the formal Fourier transformation $F^\lambda_{z,w}$ is defined as
$F^\lambda_{z,w}a(z,w)={\rm Res}_z e^{\lambda(z-w)}a(z,w)$.
Suppose $(g, F)$ is a formal distributions Lie algebra. \mbox{Define} $[a(z)\,{}_\lambda\, b(z)]=F^\lambda_{z,w}[a(z),b(w)]$ for any $a(z),b(z)\in F$. One can easily check  that this definition of \mbox{$\lambda$-brackets} satisfies (\ref{conformal}).
Given a formal distributions Lie algebra $(g, F)$, we may always include $F$ in the minimal family $F^c$ of pairwise local distributions which is closed under
the derivative $\partial$ and the $\lambda$-brackets.
Then $F^c$ is actually a Lie conformal algebra, referred to as
a {\it Lie conformal algebra of $g$}.
\begin{defi}\label{D2}\rm
A conformal module $M$ over a Lie conformal algebra $A$ is a $\C[\partial]$-module endowed with a $\lambda$-action $A\otimes M\rightarrow M[\lambda]$ such that
\begin{equation*}
\aligned
&(\partial a)\,{}_\lambda\, v=-\lambda a\,{}_\lambda\, v,\ \ \ \ \ a{}\,{}_\lambda\, (\partial v)=(\partial+\lambda)a\,{}_\lambda\, v;\\
&a\,{}_\lambda\, (b{}\,_\mu\, v)-b\,{}_\mu\,(a\,{}_\lambda\, v)=[a\,{}_\lambda\, b]\,{}_{\lambda+\mu}\, v.
\endaligned
\end{equation*}
\end{defi}
\begin{defi}\rm
A Lie conformal algebra $A$ is {\it $\Z$-graded} if $A=\oplus_{i\in \Z}A_i$, where each $A_i$ is a $\C[\partial]$-submodule
and $[A_i\,{}_\lambda\, A_j]\subset A_{i+j}[\lambda]$ for any $i,j\in \Z$.
Similarly, a conformal module $V$ over $A$ is  {\it $\Z$-graded} if $V=\oplus_{i\in \Z}V_i$, where each $V_i$ is a $\C[\partial]$-submodule and $A_i\,{}_\lambda\, V_j\subset V_{i+j}[\lambda]$ for any $i,j\in \Z$. In addition, if each $V_i$ is freely generated by one element $v_i\in V_i$ over $\C[\partial]$, we call $V$ a {\it $\Z$-graded free intermediate series module}.
\end{defi}
\section{The Lie conformal algebra $CHV$}
In this section, we  start with the Lie algebra $HV$ to construct the Lie conformal algebra $CHV$ via formal distribution Lie algebra.
Let $L_i(z)=\sum_{\a\in\Z} L_{\a,i}z^{-\a-2}$ and  $H_j(z)=\sum_{\b\in\Z} H_{\b,j}z^{-\b-1}$  for any $i,j\in\Z$. Let  $F=\{L_i(z),\,H_j(z)\,|\,i,j\in\Z\}$ be the set of $HV$-valued formal distributions.
Then we have the following result.
\begin{prop}\label{ll}We have
\begin{eqnarray*}&&
[L_i(z),L_j(w)]=(\partial_w L_{i+j}(w))\delta(z,w)+2L_{i+j}(w)\partial_w\delta(z,w),
\\
&&[L_i(z),H_j(w)]=(\partial_w H_{i+j}(w))\delta(z,w)+H_{i+j}(w)\partial_w\delta(z,w),
\\&&[H_i(z),L_j(w)]=H_{i+j}(w)\partial_w\delta(z,w),
\\&&[H_i(z),H_j(w)]=0.
\end{eqnarray*}
\end{prop}
\ni\ni{\it Proof.}\ \ Using equation (1.1), we obtain
\begin{eqnarray*}
\!\!\!\!\!\!\!\![L_i(z),L_j(w)]&\!\!\!=\!\!\!&\Big[\mbox{$\sum\limits_{\a\in\Z}$} L_{\a,i}z^{-\a-2},\mbox{$\sum\limits_{\b\in\Z}$} L_{\b,j}w^{-\b-2}\Big]\\
&\!\!\!=\!\!\!&\mbox{$\sum\limits_{\a,\b\in\Z}$}[L_{\a,i}, L_{\b,j}]z^{-\a-2}w^{-\b-2} \\
&\!\!\!=\!\!\!&\mbox{$\sum\limits_{\a,\b\in\Z}$}(\a-\b)L_{\a+\b,i+j}z^{-\a-2}w^{-\b-2} \\
&\!\!\!=\!\!\!&-\mbox{$\sum\limits_{\a,\b\in\Z}$}(\a\!+\!\b\!+\!2)L_{\a+\b,i+j}w^{-(\a+\b)-3}z^{-(\a+1)-1}w^{\a+1}\\
&\!\!\!\!\!&+\mbox{$\sum\limits_{\a,\b\in\Z}$}2(\a\!+\!1)L_{\a+\b,i+j}w^{-(\a+\b)-2}z^{-\a-2}w^{\a} \\
&\!\!\!=\!\!\!&(\partial_w L_{i+j}(w))\delta(z,w)+2L_{i+j}(w)\partial_w\delta(z,w).
\end{eqnarray*}
We also have, 
\begin{eqnarray*}
\!\!\!\!\!\!\!\![L_i(z),H_j(w)]&\!\!\!=\!\!\!&\mbox{$\sum\limits_{\a,\b\in\Z}$}[L_{\a,i}, H_{\b,j}]z^{-\a-2}w^{-\b-1}\\
&\!\!\!=\!\!\!&\mbox{$\sum\limits_{\a,\b\in\Z}$}-\b H_{\a+\b,i+j}z^{-\a-2}w^{-\b-1} \\
\end{eqnarray*}
\begin{eqnarray*}
\phantom{\!\!\!\!\!\!\!\![L_i(z),L_j(w)}
&\!\!\!=\!\!\!&\mbox{$\sum\limits_{\a,\b\in\Z}$}-(\a+\b+1)H_{\a+\b,i+j}w^{-(\a+\b)-2}z^{-(\a+1)-1}w^{\a+1}\\
&\!\!\!\!\!\!&+\mbox{$\sum\limits_{\a,\b\in\Z}$}(\a+1)H_{\a+\b,i+j}w^{-(\a+\b)-1}z^{-\a-2}w^{\a} \\
&\!\!\!=\!\!\!&(\partial_w H_{i+j}(w))\delta(z,w)+H_{i+j}(w)\partial_w\delta(z,w).
\end{eqnarray*}
The other equations follow similarly.\QED\vskip5pt

By Proposition \ref{local}, we know that $[L_i(z),L_j(w)], [L_i(z),H_j(w)], [H_i(z),L_j(w)]$ and $[H_i(z),H_j(w)]$ are local for any $i,j\in\Z$, which suggests that $F$ is a local family of formal distributions. Since the coefficients of $F$ is a basis of $HV$, we conclude that $F$ is a formal distribution Lie algebra of $HV$.

\begin{prop}\label{llL} We have
\begin{equation}\label{3.3}
\aligned
&[L_i\,{}_\lambda\, L_j]=(\partial+2\lambda)L_{i+j},
\ \ 
[L_i\,{}_\lambda\, H_j]=(\partial+\lambda)H_{i+j},
\ \ 
[H_i\,{}_\lambda\, L_j]=\lambda H_{i+j},
\ \ 
[H_i\,{}_\lambda\, H_j]=0.
\endaligned
\end{equation}
\end{prop}
\ni\ni{\it Proof.}\ \ Using Proposition \ref{ll}, we obtain
\begin{eqnarray*}
\!\!\!\!\!\!\!\![L_i(w)\,{}_\lambda\, L_j(w)]&\!\!\!=\!\!\!&
{\rm Res}_z e^{\lambda(z-w)}[L_i(z),L_j(w)]\\
&\!\!\!=\!\!\!&{\rm Res}_z e^{\lambda(z-w)}\Big((\partial_w L_{i+j}(w))\delta(z,w)+2L_{i+j}(w)\partial_w\delta(z,w)\Big)\\
&\!\!\!=\!\!\!&\partial_w L_{i+j}(w){\rm Res}_z\mbox{$\sum\limits_{j\in\Z^+}$}\lambda^k\frac{(z-w)^k\delta(z,w)}{k!}
\\&\!\!\!\!\!&+2L_{i+j}(w){\rm Res}_z\mbox{$\sum\limits_{j\in\Z^+}$}\lambda^k\frac{(z-w)^k\partial_w\delta(z,w)}{k!}\\
&\!\!\!=\!\!\!&\partial_w L_{i+j}(w){\rm Res}_z\delta(z,w)+2L_{i+j}(w){\rm Res}_z\Big(\partial_w\delta(z,w)+\lambda(z-w)\partial_w\delta(z,w)\Big)
\\
&\!\!\!=\!\!\!&
\partial_w L_{i+j}(w)+2\lambda L_{i+j}(w)
=(\partial+2\lambda)L_{i+j}(w).
\end{eqnarray*}
We also have
\begin{eqnarray*}
\!\!\!\!\!\!\!\![L_i(w)\,{}_\lambda\, H_j(w)]&\!\!\!=\!\!\!&
{\rm Res}_z e^{\lambda(z-w)}[L_i(z),H_j(w)]\\
&\!\!\!=\!\!\!&{\rm Res}_z e^{\lambda(z-w)}\Big((\partial_w H_{i+j}(w))\delta(z,w)+H_{i+j}(w)\partial_w\delta(z,w)\Big)\\
&\!\!\!=\!\!\!&-\partial_w H_{i+j}(w){\rm Res}_z\mbox{$\sum\limits_{j\in\Z^+}$}\lambda^k\frac{(z-w)^k\delta(z,w)}{k!}
\\&\!\!\!\!\!&+H_{i+j}(w){\rm Res}_z\mbox{$\sum\limits_{j\in\Z^+}$}\lambda^k\frac{(z-w)^k\partial_w\delta(z,w)}{k!}\\
&\!\!\!=\!\!\!&\partial_w H_{i+j}(w){\rm Res}_z\delta(z,w)+L_{i+j}(w){\rm Res}_z\Big(\partial_w\delta(z,w)+\lambda(z-w)\partial_w\delta(z,w)\Big)\\
&\!\!\!=\!\!\!&\partial_w H_{i+j}(w)+\lambda H_{i+j}(w)
=(\partial+\lambda)H_{i+j}(w).
\end{eqnarray*}
Using the same method, we  
 obtain 
\begin{equation*}
[H_i(w)\,{}_\lambda\, L_j(w)]=\lambda H_{i+j}(w),
[H_i(w)\,{}_\lambda\, H_j(w)]=0.
\eqno\mbox{\QED}\end{equation*}

\begin{prop}\label{lfom1}
Let $CHV$ be a free $C[\partial]$-module with basis $\C[\partial]$-basis $\{L_i, H_i\,|\,i\in\Z\}$. Then $CHW$ is a Lie conformal algebra with $\lambda$-brackets defined as in \eqref{3.3}.
\end{prop}

Note that $CHV$ is a $\Z$-graded Lie conformal algebra in the sense $CHV=\oplus_{i\in\Z} \ (CHV)_{i}$, where
$(CHV)_{i}=C[\partial]{L_i}\oplus \C[\partial]{H_i}$.

\section{Conformal derivations of $CHV$}
Suppose $A$ is a Lie conformal algebra.
A linear map $\phi_\lambda: A\rightarrow A[\lambda]$ is called a conformal derivation if the following equalities hold:
\begin{equation}\label{4.1}
\aligned
&\phi_\lambda(\partial v)=(\partial+\lambda)\phi_\lambda(v),\ \ \ \
&\phi_\lambda([a\,{}_\mu \,b])=[(\phi_\lambda a)\,{}_{\lambda+\mu} \,b]+[a\,{}_\mu \,(\phi_\lambda b)].
\endaligned
\end{equation}
We often write $\phi$ instead of $\phi_\lambda$ for simplicity.

It can be easily verified that for any $x\in A$, the map ${\rm ad}_x$, defined by $({\rm ad}_x)_\lambda y= [x\, {}_\lambda\, y]$ for   $y\in A$, is a conformal derivation of $A$.
All conformal derivations of this kind are called {\it inner}.
Denote by ${\rm CDer\,}(CHV)$
and ${\rm CInn\,}(CHV)$ the vector spaces of all conformal derivations and inner conformal derivations of $CHV$, respectively.
Assume $D\in {\rm CDer\,} (CHV)$.
Define $D^i(L_j)=\pi_{i+j} D(L_j),\,D^i(H_j)=\pi_{i+j} D(H_j)$ for any $j\in\Z$, where in general $\pi_{i}$ is the natural projection from $$\C[\lambda]\otimes CHV\cong \OP{k\in\Z}\C[\partial,\lambda]L_k\oplus\OP{j\in\Z}\C[\partial,\lambda]H_j,$$ onto $\C[\partial,\lambda]{L_{i}}\oplus\C[\partial,\lambda]{H_{i}}$.
Then $D^i$ is a conformal derivation and $D=\sum_{i\in\Z} D^i$ in the sense that for any $x\in CHV$ only finitely many $D^i_\lambda(x)\neq0$.
Let ${({\rm CDer\,}(CHV))}^c$ be the space of conformal derivations of degree $c$, i.e.,
 $${({\rm CDer\,}(CHV))}^c=\{D\in {\rm CDer\,}(CHV)\,|\, D_\lambda((CHV)_i)\subset (CHV)_{i+c}[\lambda]\}.$$
For $a\in\C, c\in\Z$, we define
$(D^c_a) {}_\lambda(L_i)=aH_{i+c}$ and $ (D^c_a) {}_\lambda(H_i)=0$ for any $i\in\Z$.
Then \linebreak
$D^c_a\in {({\rm CDer\,}(CHV))}^c$.

\begin{lemm}\label{lc1}
Assume $D^c\in {({\rm CDer\,}(CHV))}^c$. There exists an $a\in \C$ such that $D^c-D^c_a\in {\rm CInn\,}(CHV)$.
\end{lemm}
\ni\ni{\it Proof.}\ \   Assume $D^c_\lambda(L_i)=f_{1,i}(\partial,\lambda)L_{i+c}+f_{2,i}(\partial,\lambda)H_{i+c}$.
Since $[L_0\ {}_\mu \ L_i]=(\partial+2\mu)L_i$,
one has
\begin{equation}\label{c1}
(\partial+\lambda+2\mu)f_{1,i}(\partial,\lambda)=(\partial+2\lambda+2\mu)f_{1,0}(-\lambda-\mu,\lambda)+(\partial+2\mu)f_{1,i}(\partial+\mu,\lambda),
\end{equation}
and
\begin{equation}\label{c2}
(\partial+\lambda+2\mu)f_{2,i}(\partial,\lambda)=(\lambda+\mu)f_{2,0}(-\lambda-\mu,\lambda)+(\partial+\mu)f_{2,i}(\partial+\mu,\lambda).
\end{equation}
Setting $\mu=0$ in (\ref{c1}) and (\ref{c2}), one gets
\begin{equation}\label{4.2}
\lambda f_{1,i}(\partial,\lambda)=(\partial+2\lambda)f_{1,0}(-\lambda,\lambda),
\ \ \ \ \ 
\lambda f_{2,i}(\partial,\lambda)=\lambda f_{2,0}(-\lambda,\lambda).
\end{equation}
From the first equation of (\ref{4.2}), $\lambda$ is a factor of  $f_{1,0}(-\lambda,\lambda)$ in the unique factorization ring $C[\partial,\lambda]$. Setting $g_1(\lambda)=\frac{f_{1,0}(\lambda,-\lambda)}{\lambda}$ and replacing $D^c$ by $D^c_\lambda-{\rm ad}_{-{g_1(\partial)}L_c}$,
we can suppose $D^c_\lambda(L_i)=f_{2,i}(\partial,\lambda)H_{i+c}$. Then by the second equation of (\ref{4.2}), there exists a polynomial
$g(\lambda)\in\C[\lambda]$ such that $f_{2,i}(\partial,\lambda)=g(\lambda)$ for  $i\in \Z$.
Assume $g(\lambda)=\lambda h(\lambda)+a$. Then $[h(-\partial)H_c\, {}_\lambda \,L_i]=\lambda h(\lambda)H_{i+c}$.
Replacing $D^c$ by $D^c-{\rm ad}_{{h(-\partial)}H_c}$, we can assume $D^c_\lambda(L_i)=aH_{i+c}$ for any $i\in \Z$.

Assume $D^c_\lambda(H_i)=f_{3,i}(\partial,\lambda)L_{i+c}+f_{4,i}(\partial,\lambda)H_{i+c}$. Let $D^c_\lambda$ act on $[L_0\ {}_\mu \ H_i]=(\partial+\mu)H_i$,
one has
\begin{equation*}
(\partial+2\mu)f_{3,i}(\partial+\mu,\lambda)=(\partial+\mu+\lambda)f_{3,i}(\partial,\lambda),\ \
%
(\partial+\mu)f_{4,i}(\partial+\mu,\lambda)=(\partial+\mu+\lambda)f_{4,i}(\partial,\lambda).
\end{equation*}
Setting $\mu=0$, 
one gets
$
\ f_{3,i}(\partial,\lambda)=f_{4,i}(\partial,\lambda)=0.
$ 
Hence  $D^c_\lambda(H_i)=0$.\QED
\begin{lemm}\label{lc2}We have
${\rm CDer\,}(CHV)=\oplus_{c\in\Z}{({\rm CDer\,}(CHV))}^c$.
\end{lemm}
\ni\ni{\it Proof.}\ \  From Lemma \ref{lc1}, we have $D=\sum_{c\in\Z} D^c$, where $D^c=ad_{g_c(\partial)L_c+h_c(\partial)H_c}+D^c_{a_c}$ for some $g_c(\partial), h_c(\partial)\in\C[\partial]$ and $a_c\in\C$.
Since
\begin{equation*}
D^c{}_\lambda(L_0)=[(g_c(\partial)L_c+h_c(\partial)H_c)\, {}_\lambda\, L_0]+D^c_{a_c}{}_\lambda(L_0)=(\partial+2\lambda)g_c(-\lambda)L_c+\lambda h_c(-\lambda)H_c+a_c H_c,
\end{equation*}
we see that $D^c{}_\lambda(L_0)=0$ implies $D^c=0$. Consequently, we deduce that $D=\sum_{c\in\Z} D^c$ is a finite sum.\QED\vskip5pt

Denote $$\C^\infty=\{\vec{a}=(a_c)_{c\in\Z}\,|\,a_c\in\C\mbox{ and }a_c=0\ \mbox{for\ all\ but\ finitely\ many}\ c's\}.$$
For each $\vec{a}\in \C^\infty$, we define $D_{\vec{a}}\, {}_\lambda\, (L_i)=\sum a_c H_{i+c}$ and $ D_{\vec{a}}\, {}_\lambda\,(H_i)=0$ for all $i\in\Z$.
Then $D_{\vec{a}}\in {\rm CDer\,}(CHV)$. It  can be easily verified that $D_{\vec{a}}\in {\rm CInn\,}(CHV)$ implies $D_{\vec{a}}=0$.
We also denote by $\C^\infty$ the space of such conformal derivations.
Therefore, the following theorem follows from Lemma \ref{lc1} and Lemma \ref{lc2}.
\begin{theo}We have
${\rm CDer\,}(CHV)={\rm CInn\,}(CHV)\oplus \C^\infty$.
\end{theo}
\section{Rank one conformal modules over $CHV$}
Suppose  $M$ is a free conformal module of rank one over $CHV$.
We may write $M=\C[\partial]v$ and assume $L_i\,{}_\lambda\, v=f_i(\partial,\lambda)v$,  $H_j\,{}_\lambda\, v=g_j(\partial,\lambda)v$, where $f_i(\partial,\lambda),g_j(\partial,\lambda)\in\C[\partial,\lambda]$. We will compute the coefficients $f_i(\partial,\lambda),g_j(\partial,\lambda)$ in the rest of this section.

For the Virasoro conformal algebra $CVir$, it is well known that all the free nontrivial $CVir$-modules of rank one over $\C[\partial]$ are the following ones $(a,b\in\C)$,
\begin{equation}\label{m0}
M^\prime_{a,b}=\C[\partial]v,\ \ \ \ L{}_\lambda \,v=(\partial+a\lambda+b)v.
\end{equation}
The module $M^\prime_{a,b}$ is irreducible if and only if $a\neq 0$. The module $M^\prime_{0,b}$ contains a unique nontrivial submodule $(\partial+b)M^\prime_{0,b}$  isomorphic to $M^\prime_{1,b}$. It was proved that the modules $M^\prime_{a,b}$ with $a\neq 0$ exhaust all finite irreducible nontrivial $CVir$-modules in \cite{CK}.

First, by 
\cite{WCY}, we have the following result.
\begin{lemm}\label{lfo3}
There exist $a,b,c\in\C$ such that 
$
f_i(\partial,\lambda)=c^i(\partial+a\lambda+b)$ 
for  $i\in\Z$.\end{lemm}

We continue to compute $g_i(\partial,\lambda)$ in the following Lemma.
\begin{lemm}\label{m20}
There exists $d\in\C$ such that
$
g_i(\partial,\lambda)=dc^i.
$
\end{lemm}
\ni\ni{\it Proof.}\ \
 A direct computation shows that
\begin{eqnarray*}&&
[H_i\,{}_\lambda\, H_j]\,{}_{\lambda+\mu} v=H_i\,{}_\lambda\,(H_j\, {}_\mu v)-H_j\,{}_\mu\, (H_i \,{}_\lambda v)=0,
\\&& 
H_i\,{}_\lambda\,(H_j\, {}_\mu v)=H_i\,{}_\lambda (g_j(\partial,\mu)v)=g_j(\partial+\lambda,\mu)H_i {}_\lambda v=g_j(\partial+\lambda,\mu)g_i(\partial,\lambda)v,
\end{eqnarray*}
and
\begin{equation*}
H_j\,{}_\mu\, (H_i \,{}_\lambda v)=g_i(\partial+\mu,\lambda)g_j(\partial,\mu)v.
\end{equation*}
Thus 
$
g_j(\partial+\lambda,\mu)g_i(\partial,\lambda)=g_i(\partial+\mu,\lambda)g_j(\partial,\mu).
$ 
Comparing the coefficients of $\lambda$, we obtain that
$g_i(\partial,\lambda)$ is independent of the variable $\partial$ and denote
$
g_i(\lambda)=g_i(\partial,\lambda)\mbox{ \ for \ }i\in\Z.
$ 
Using 
\begin{eqnarray*}
\!\!\!\!\!\!\!\!\!\!\!\!\!\!
&&[L_i\,{}_\lambda\, H_j]\,{}_{\lambda+\mu} v=((\partial+\lambda) H_{i+j}){}_{\lambda+\mu} v=-\mu g_{i+j}(\lambda+\mu)v,
\\
\!\!\!\!\!\!\!\!\!\!\!\!\!\!&&[L_i\,{}_\lambda\, H_j]\,{}_{\lambda+\mu} v=L_i\,{}_\lambda\,(H_j\, {}_\mu v)-H_j\,{}_\mu\, (L_i \,{}_\lambda v)
=g_{i+j}(\mu)(f_i(\partial,\lambda )-f_i(\partial+\mu,\lambda))v=-\mu c^i g_{j}(\mu)v.
\end{eqnarray*}
we obtain 
\begin{equation*}
\mu g_{i+j}(\lambda+\mu)=\mu c^i g_{j}(\mu),\ \
\ \ \ 
g_i(\lambda)=d_i.
\end{equation*}
Thus $
d_{i+j}=c^i d_j.
$ 
Denote $d=d_0$.
Then $
d_{i}=dc^i.
$ 
Thus $
H_i\, {}_\lambda \,v=g_i(\lambda)v=dc^i v.
$ 
This implies the result.\QED\vskip5pt
Now we get the main result of this section.
\begin{theo}\label{theo1}
A  nontrivial free conformal module of rank one over $CHV$ is isomorphic to $M_{a,b,c,d}$ for some $a,b,d\in\C, c\in\C^*$, where $M_{a,b,c,d}=\C[\partial]v$ and $\lambda$-actions are given by
\begin{equation}
L_i\, {}_\lambda \,v=c^i(\partial+a\lambda+b)v,\ \ H_i\, {}_\lambda \,v=dc^i v.
\end{equation}
Furthermore, $M_{a,b,c,d}$ is irreducible if and only if $a\neq0$.
\end{theo}

\section{The $\Z$-graded free intermediate series modules over $CHV$}

In this section, we give a classification of $\Z$-graded free intermediate series modules over $CHV$.
Let $V$ be an arbitrary $\Z$-graded free intermediate series module over $CHV$. Then $V=\oplus_{i\in\Z}V_i$, where each $V_i$
is freely generated by some element $v_i\in V_i$ over $\C[\partial]$.
For any $i,j\in\Z$, we denote \begin{equation*}
L_i\,{}_\lambda\, v_k=f_{i,k}(\partial,\lambda)v_{i+k}, \ \ \ H_j\,{}_\lambda\, v_k=g_{j,k}(\partial,\lambda)v_{j+k}\mbox{ for some }f_{i,k}(\partial,\lambda),g_{j,k}(\partial,\lambda)\in\C[\partial,\lambda].\end{equation*}
We call $\{f_{i,k}(\partial,\lambda),\,g_{j,k}(\partial,\lambda)\,|\,i,j,k\in\Z\}$ the {\it structure
coefficients} of $V$ related to the $\C[\partial]$-basis $\{v_i\,|\,i\in\Z\}$. Then the conformal module structure on $V$ is determined
if and only if all of its structure coefficients are specified.

Since $V$ is also a $\Z$-graded free intermediate series module over $CW$ (c.f. \eqref{CWAAAA}),
we can employ a result of 
\cite{WCY}. 
First we introduce two classes of $\Z$-graded free intermediate series module over $CW$:
Given $a,b\in\C$, let $V^\prime_{a,b}=\oplus_{i\in\Z}\C[\partial]v_i$ and
define
\begin{equation}\label{act1111}
L_i\,{}_\lambda\, v_j=(\partial+a\lambda+b)v_{i+j}.
\end{equation}
Then $V^\prime_{a,b}$ is a $\Z$-graded free intermediate series module over $CW$.

We denote $\{0,1\}^\infty$ to be the set of sequences
$A=\{a_i\}_{i\in\Z}$  with $a_i\in\{0,1\}$ for any $i\in\Z$. Let $A\in\{0,1\}^\infty$ and $b\in\C$.
We can construct a nontrivial $\Z$-graded free intermediate series module $V^\prime_{A,b}$, where $V^\prime_{A,b}=\oplus_{i\in\Z}\C[\partial]v_i$ and
the $\lambda$-actions are given by
\begin{equation}\label{act11112}
L_i\,{}_{{}_\lambda} v_j=
\begin{cases}
(\partial+b)v_{i+j} &\ \mbox{if}\  (a_j,a_{i+j})=(0,0);\\[4pt]
(\partial+b+\lambda)v_{i+j}&\  \mbox{if} \ (a_j,a_{i+j})=(1,1);\\[4pt]
v_{i+j} &\ \mbox{if} \ (a_j,a_{i+j})=(0,1);\\[4pt]
(\partial+b)(\partial+b+\lambda)v_{i+j}&\  \mbox{if} \ (a_j,a_{i+j})=(1,0).
\end{cases}
\end{equation}

By 
\cite{WCY}, we have the following result.
\begin{lemm}
Assume $V^\prime$ is a nontrivial $\Z$-graded free intermediate series module over $CW$.
Then $V^\prime$ is isomorphic to $V^\prime_{a,b}$ defined by \eqref{act1111} or $V^\prime_{A,b}$ defined by \eqref{act11112}.\QED
\end{lemm}

From Definition \ref{D2}, one can obtain the following result.
\begin{lemm}\label{lf1}
The structure coefficients $f_{i,j}(\partial,\lambda)$ and $g_{i,j}(\partial,\lambda)$ of
$V$ satisfy the following:
\begin{eqnarray}\label{f148}&&
\mu g_{i+j,k}(\partial,\lambda+\mu)=f_{i,k}(\partial+\mu,\lambda)g_{j,i+k}(\partial,\mu)-g_{j,k}(\partial+\lambda,\mu)f_{i,j+k}(\partial,\lambda),
\\ 
\label{f149}
&& g_{j,k}(\partial+\lambda,\mu)g_{i,j+k}(\partial,\lambda)=g_{i,k}(\partial+\mu,\lambda)g_{j,i+k}(\partial,\mu).
\end{eqnarray}
\end{lemm}
\ni\ni{\it Proof.}\ \   One can easily check  that
\begin{equation*}
[L_i\,{}_\lambda\, H_j]\,{}_{\lambda+\mu} \,v_k=(\partial H_{i+j}+\lambda H_{i+j})\,{}_{\lambda+\mu} \,v_k=-\mu g_{i+j,k}(\partial,\lambda+\mu)v_{i+j+k}.
\end{equation*}
Similarly we  have the following two equations:
\begin{equation*}
\aligned
L_i\,{}_\lambda \,(H_j\, {}_\mu \,v_k)&=L_i\,{}_\lambda \,(g_{j,k}(\partial,\mu)v_{j+k})
=g_{j,k}(\partial+\lambda,\mu)f_{i,j+k}(\partial,\lambda)v_{i+j+k};\\
H_j\,{}_\mu \,(L_i \,{}_\lambda\, v_k)&=f_{i,k}(\partial+\mu,\lambda)g_{j,i+k}(\partial,\mu)v_{i+j+k}.
\endaligned
\end{equation*}
 One can also check that
\begin{equation*}
[H_i\,{}_\lambda\, H_j]\,{}_{\lambda+\mu} \,v_k=H_i\,{}_\lambda \,(H_j\, {}_\mu \,v_k)-H_j\,{}_\mu \,(H_i \,{}_\lambda\, v_k)=0.
\end{equation*}
Then 
we  have the following two equations:
\begin{equation*}
\aligned
H_i\,{}_\lambda \,(H_j\, {}_\mu \,v_k)&=H_i\,{}_\lambda \,(g_{j,k}(\partial,\mu)v_{j+k})
=g_{j,k}(\partial+\lambda,\mu)g_{i,j+k}(\partial,\lambda)v_{i+j+k};\\
H_j\,{}_\mu \,(H_i \,{}_\lambda\, v_k)&=g_{i,k}(\partial+\mu,\lambda)g_{j,i+k}(\partial,\mu)v_{i+j+k}.
\endaligned
\end{equation*}
Now the result 
follows from the defining relations of a conformal module.\QED
\begin{lemm}\label{lf21} The
 $g_{0,k}(\partial,\lambda)$ is independent of the variable $\partial$.
\end{lemm}
\ni\ni{\it Proof.}\ \
Letting $i=j=0$ in (\ref{f149}),  we obtain
\begin{eqnarray*}
g_{0,k}(\partial+\lambda,\mu)g_{0,k}(\partial,\lambda)=g_{0,k}(\partial+\mu,\lambda)g_{0,k}(\partial,\mu),
\\ 
\deg_{y} g_{0,k}(x,y)=\deg_{x} g_{0,k}(x,y)+\deg_{y} g_{0,k}(x,y).
\end{eqnarray*}
Hence $\deg_{x} g_{0,k}(x,y)=0$.\QED\vskip5pt

By  equation (\ref{f148}), a nontrivial conformal module over $CHV$ is also a nontrivial conformal over $CW$.
\begin{prop}\label{P110}
If  $V\cong V^\prime_{a,b}$ as $CW$-modules for some $a,b\in\C$, then $g_{j,k}(\partial,\lambda)=c$ for some $c\in\C$ and for all $j,k\in\Z$.
In particular, $H_j\,{}_\lambda\, v_k=cv_{j+k}$ and denote this module by $V_{a,b,c}$.
\end{prop}
\ni\ni{\it Proof.}\ \
If $V=V^\prime_{a,b}$, then $f_{i,j}(\partial,\lambda)=(\partial+a\lambda+b)$.
Letting $i=j=0$ in (\ref{f148}), we  obtain
\begin{equation}\label{f122}
\mu g_{0,k}(\lambda+\mu)=f_{0,k}(\partial+\mu,\lambda)g_{0,k}(\mu)-g_{0,k}(\mu)f_{0,k}(\partial,\lambda),
\end{equation}
where we have used the previous lemma to re-denote $g_{0,k}(\partial,\lambda)$ by $g_{0,k}(\mu)$.
Thus 
$
\mu  g_{0,k}(\lambda+\mu)=\mu  g_{0,k}(\mu),$ 
Hence $ g_{0,k}(\lambda+\mu)= g_{0,k}(\mu)$. Therefore,  we can assume $g_{0,k}(\partial,\lambda)=x_{k}$ for some $x_{k}\in\C$ and for all $k\in\Z$.
Letting $j=0$ in (\ref{f148}),  we obtain
\begin{equation*}
\mu g_{i,k}(\partial,\lambda+\mu)=f_{i,k}(\partial+\mu,\lambda)g_{0,k}(\partial,\mu)-g_{0,i+k}(\partial+\lambda,\mu)f_{i,k}(\partial,\lambda),
\end{equation*}i.e.,
\begin{equation}\label{f124}
\mu g_{i,k}(\partial,\lambda+\mu)=(\partial+a\lambda+b)(x_{i+k}-x_{k})+\mu x_{i+k}.
\end{equation}
By considering coefficients of $\mu^i$ for $i=0,1$ in
\eqref{f124},
we immediately obtain
%
$x_i=c$ for some $c\in\C$, and further, $g_{i,k}(\partial,\lambda)=c$.\QED\vskip5pt

We can similarly obtain the following proposition.
\begin{prop}\label{P111}
If  $V\cong V^\prime_{A,b}$ as $CW$-modules for some $A\in\{0,1\}^\infty,\,b\in\C$,
then $g_{j,k}(\partial,\lambda)=c$ for some $c\in\C$ and for all $j,k\in\Z$.
In particular, $H_j\,{}_\lambda\, v_k=cv_{j+k}$ and denote this module by $V_{A,b,c}$. \end{prop}

Now 
we get the main result of this section.
\begin{theo}
Assume $V$ is a nontrivial $\Z$-graded free intermediate series module over $CHV$.
Then $V$ is isomorphic to either $V_{a,b,c}$ defined by \eqref{act1111} and Proposition \ref{P110} for some $a,b,c\in\C$, or else $V_{A,b,c}$ defined by \eqref{act11112} and Proposition \ref{P111} for some $A\in\{0,1\}^\infty,\,b, c\in\C$.
\end{theo}
%
\small

\end{document}